\def\Hl {{H_{2,2}}}
\def\Al {{\cA_{1,3}(2)}}
\def\bAl {{\Bar{\cA}_{1,3}(2)}}
\def\Ztwo {{\ZZ_2}}
\def\BKm {\mathop{\Tilde{\rm Km}}\nolimits}
\def\pL {{\phi_{|\tilde\cL|^-}}}
\def\Gl {{\Gamma_{1,3}(2)}}
\def\Gpara {{\Gamma_{1,3}}}
\def\SCD {{\Sigma_{(C:D)}}}
\def\tS {{\tilde\Sigma_{(C:D)}}}

\magnification=1200
\input amssym.def
\input amssym
\def\AA {{\Bbb A}}
\def\CC {{\Bbb C}}

\def\HH {{\Bbb H}}

\def\PP {{\Bbb P}}
\def\QQ {{\Bbb Q}}
\def\RR {{\Bbb R}}
\def\SS {{\Bbb S}}
\def\ZZ {{\Bbb Z}}
\def\contin {\subseteq}
\def\Tilde{\widetilde}
\def\Hat{\widehat}
\def\Bar{\overline}

\def\cross {\times}
\def\sqtimes {\boxtimes}

\def\NS {\mathop{\rm NS}\nolimits}
\def\Km {\mathop{\rm Km}\nolimits}

\def\Aut {\mathop{\rm Aut}\nolimits}
\def\Hom {\mathop{\rm Hom}\nolimits}

\def\Sp {\mathop{\rm Sp}\nolimits}
\def\SL {\mathop{\rm SL}\nolimits}

\def\Pic {\mathop{\rm Pic}\nolimits}
\def\Jac {\mathop{\rm Jac}\nolimits}
\def\Grass {\mathop{\rm Grass}\nolimits}

\def\im {\mathop{\rm Im}\nolimits}

\def\rk {\mathop{\rm rk}\nolimits}
\def\cA {{\cal A}}
\def\cF {{\cal F}}
\def\cG {{\cal G}}
\def\cH {{\cal H}}
\def\cL {{\cal L}}
\def\cM {{\cal M}}

\def\cO {{\cal O}}

\def\bdu {{\bf u}}

\def\bdx {{\bf x}}

\def\To{\longrightarrow}

\def\Sum{\sum\limits}

\def\tens{\otimes}

\def\pf{\noindent{\sl Proof: }}

\def \hf {{{1}\over{2}}}
\def\rat{\mathrel{{\hbox{\kern2pt\vrule height2.45pt depth-2.15pt
 width2pt}\kern1pt {\vrule height2.45pt depth-2.15pt width2pt}
  \kern1pt{\vrule height2.45pt depth-2.15pt width1.7pt\kern-1.7pt}
   {\raise1.4pt\hbox{$\scriptscriptstyle\succ$}}\kern1pt}}}
\def\qed{\vrule width5pt height5pt depth0pt\par\smallskip}
\def\surj{\to\kern-8pt\to}
\outer\def\startsection#1\par{\vskip0pt
 plus.3\vsize\penalty-10\vskip0pt
  plus-.3\vsize\bigskip\vskip\parskip\message{#1}
   \leftline{\bf#1}\nobreak\smallskip\noindent}
\def\imic{\cong}
\outer\def\thm #1 #2\par{\medbreak
  \noindent{\bf Theorem~#1.\enspace}{\sl#2}\par
   \ifdim\lastskip<\medskipamount \removelastskip\penalty55\medskip\fi}
\outer\def\prop #1 #2\par{\medbreak
  \noindent{\bf Proposition~#1.\enspace}{\sl#2}\par
   \ifdim\lastskip<\medskipamount \removelastskip\penalty55\medskip\fi}
\outer\def\lemma #1 #2\par{\medbreak
  \noindent{\bf Lemma~#1.\enspace}{\sl#2}\par
   \ifdim\lastskip<\medskipamount \removelastskip\penalty55\medskip\fi}
\outer\def\corollary #1 #2\par{\medbreak
  \noindent{\bf Corollary~#1.\enspace}{\sl#2}\par
   \ifdim\lastskip<\medskipamount \removelastskip\penalty55\medskip\fi}
\def\deep #1 {_{\lower5pt\hbox{$#1$}}}
\font\headfont=cmb10 scaled\magstep2
\font\smallfont=cmr8
\baselineskip=15pt
\raggedbottom

\centerline{\headfont Heisenberg-invariant Kummer surfaces}
\bigskip
\centerline{K.~Hulek, I.~Nieto \& G.K.~Sankaran}

{\narrower\medskip\noindent {\smallfont We study, from the point of view of
abelian and Kummer surfaces and their moduli, the special quintic threefold
known as Nieto's quintic. It is in some ways analogous to the Segre cubic
and the Burkhardt quartic and can be interpreted as a moduli space of
certain Kummer surfaces. It contains $30$ planes and has $10$ singular
points: we describe how some of these arise from bielliptic and product
abelian surfaces and their Kummer surfaces.\smallskip}}

\medskip
\noindent
In this paper we study, from the point of view of abelian surfaces and
their moduli, the special quintic threefold $N$ first described in [Ni1]
and [BN], known as Nieto's quintic\footnote {${}^1$}{{\smallfont This
terminology is used in [Hun] and elsewhere and is adopted here by a
majority vote of the authors.}}. Nieto's quintic is in some ways analogous
to the Segre cubic and the Burkhardt quartic and has a rich but still
largely unexplored geometry. It has a double cover $\tilde N$ which is
birationally equivalent to the moduli space $\Al$ of $(1,3)$-polarised
abelian surfaces with a level-$2$ structure. $N$~contains $30$ planes and
has $10$ singular points: we aim to understand how these arise from the
abelian surfaces and their Kummer surfaces. It turns out that $15$ of the
planes are related to degenerate abelian surfaces and this aspect is
studied in the companion paper~[HNS]. The other $15$ planes come from a
certain (reducible) Humbert surface in~$\Al$. The $10$ singular points
arise similarly from another such Humbert surface. These two Humbert
surfaces correspond to precisely those abelian surfaces whose minimally
resolved Kummer surface is not embedded by the anti-invariant part of twice
the polarisation.
\bigskip

\vbox{
\baselineskip=10pt
\noindent{\smallfont Acknowledgments: The orginal
impetus for this work came from discussions between GKS and IN while
GKS was visiting Mexico with support from the EU International
Scientific Cooperation Programme Contract No. CI1*-CT93-0031. Most of
the work on this paper was done during the visit of IN to Bath,
supported by EPSRC grant number GR/L27534, and Hannover in the autumn
of 1996. The second author was also supported by Conacyt project
0325P-E ``Moduli of Polarized Abelian Varieties''. KH and GKS also
acknowledge support from the EU HCM network AGE.}}

\noindent
1991 Mathematics Subject Classification 14K10; 14J25, 14J30, 14N05

\startsection 1. Background and notation

We need to fix some notation for the projective varieties discussed in
[BN] and we also need to make careful definitions concerning moduli of
abelian surfaces.

$N$ is most conveniently defined as the subvariety of $\PP^5$ given by
$$
\sum_{i=0}^5u_i=\sum_{i=0}^5u_i^{-1}=0
$$
where the $u_i$ are homogeneous coordinates on $\PP^5$. As is shown in
[BN], $N$ is singular at the ten points which are equivalent
to $(1:1:1:-1:-1:-1)$ under the permutation action of the symmetric
group $\SS_6$, and along the twenty lines $L_{ijk}$ given by
$u_i=u_j=u_k=0$ which we will call desmic lines or D-lines. 

$N$ contains $30$~planes forming two $\SS_6$-orbits of size 15: the
planes
$$
u_{\sigma(0)}+u_{\sigma(1)}=u_{\sigma(2)}+u_{\sigma(3)}
=u_{\sigma(4)}+u_{\sigma(5)}=0, \qquad\sigma\in\SS_6
$$
which we call the S-planes, and the planes $F_{ij}$ given by
$$
u_i=u_j=0
$$
which we call the V-planes.

The double cover $\tilde\PP^5\to\PP^5$ branched along the coordinate
hyperplanes induces by pullback a double cover $\nu:\tilde N\to N$ which
is branched only along the V-planes.

We recall the definition of the Heisenberg group $\Hl$. It is a
central extension
$$
0\To\mu_2\To\Hl\To\ZZ_2^{\oplus 4}\To 0,
$$
where $\mu_2$ is the group of square roots of unity, which acts on
$V_{2,2}=\Hom(\Ztwo\oplus\Ztwo,\CC)$ via the usual Schr\"odinger
representation.

This induces an action of $\Hl$ on $\PP^3$. As explained in [BN],
$N$~is the closure of the locus in $\PP^4=\{\Sum u_i =0\}$ that
parametrises smooth $\Hl$-invariant quartic surfaces $X\contin\PP^3$ that
contain a line. The S- and V-planes correspond to singular quartic
surfaces~$X$. In this paper we shall be concerned mainly with the
V-planes.

We also consider the threefold $M\contin\Grass(2,4)$ which
parametrises lines in $\PP^3$ contained in some Heisenberg-invariant
quartic surface $X\contin\PP^3$. In appropriate homogeneous
coordinates $x_i$ on $\PP^5$ it is given by
$$
\sum_{i=0}^5x^2_i=\sum_{i=0}^5x_i^{-2}=0
$$
and has an action of $\Hl$, as well as a map $M\to N$ given by
$u_i=x^2_i$.

Turning to abelian surfaces, it is necessary to be careful with
definitions. What follows is all standard but is restated here for
clarity. An abelian surface is a projective algebraic group of
dimension~$2$. It is thus a projective complex torus $A$ with a group
structure, including in particular a distinguished point $0\in A$ and
an involution $\iota:A\to A$ as well as a multiplication law. A
polarisation of type $(d_1,d_2)$ on $A$, where $d_1$, $d_2$ are
positive integers and $d_1|d_2$, is an algebraic equivalence class (or
cohomology class) $H\in\NS(A)\contin H^2(A,\ZZ)$ which is the first
Chern class of an ample line bundle $\cL$ of type $(d_1,d_2)$. This
makes sense even if we do not fix an origin $0\in A$. The class
${{1}\over{d_1}}H$ is also an integral class and defines a
polarisation of type $(1,d_2/d_1)$, but $\cL^{1/d_1}$ is not well
defined: if $\cM$ is a line bundle such that $\cM^{d_1}\imic\cL$ then
$(t_\bdx^*\cM^{d_1})\imic\cL$ also for any $d_1$-torsion point
$\bdx\in A$.

An abelian torsor is a principal homogeneous space for an abelian
variety. By a symmetric torsor we mean a torsor $Y$ for an abelian
variety $A$ on which a faithful action of $\Aut A$ is also specified;
in particular, there is an involution $-1$ on~$Y$ with $16$~fixed
points, which one might refer to as possible origins (in that the
choice of one of them as origin would make $Y$ into an abelian
variety) or $2$-torsion points. An example of a symmetric torsor is
$Y=\Pic^H A$, for a general $H\in NS(A)$: if $\sigma\in\Aut A$ then it
acts on $Y$ by pullback.

The following easy fact will be useful to us.

\lemma 1.1 If $H$ is a polarisation of type $(2,2d)$ then there is a
unique symmetric line bundle $\cL_H$ with $c_1(\cL_H)=H$ such that
$\cL$ is totally symmetric, that is, the square of a symmetric line
bundle of type $(1,d)$.~\qed

Suppose we have an abelian variety $A$ and a symmetric line bundle
$\cL$ on $A$ such that the group $K(\cL)=\{\bdx\in
A|t_\bdx^*\cL\imic\cL\}$ contains the group ${}_2A$ of $2$-torsion
points of $A$. As in [Mu] and [LB], we define a theta structure of
level~$(2,2)$ on $A$ to be an isomorphism between the
Heisenberg group $\Bar H_{2,2}$, which is an extension
$$
0\To\CC^*\To\Bar H_{2,2}\To\ZZ_2^{\oplus 4}\To 0,
$$
and the group $\cG_{2,2}(\cL)$ of theta characteristics for $\cL$, given by
$$
0\To\CC^*\To\cG_{2,2}\To{}_2A\To 0,
$$
which is the identity on $\CC^*$. Such a choice induces an action of
$\Hl$ on $\cL$ via the Schr\"odinger representation of $\Hl$ and also
a choice of symplectic basis for ${}_2A$ over $\Ztwo$ with respect to
the symplectic (Weil) form induced by $H=c_1(\cL)$. In particular it
depends only on $c_1(\cL)$, not on $\cL$ itself. The choice of
symplectic basis is a level-$(2,2)$ structure in the sense of [Mu],
often called simply a level-$2$ structure. We shall sometimes abuse
notation by speaking of a level-$2$ structure for an abelian surface
with a polarisation of type $(1,3)$, when in fact $K(\cL)$ does not
contain ${}_2A$. If we do so we mean a level-$2$ structure for $A$ with
the polarisation $c_1(\cL^{\otimes 2})$, which makes sense.

We denote by $\Al$ the (coarse) moduli space of abelian surfaces with
a polarisation of type $(1,3)$ and a level-$2$ structure. It makes no
difference to the moduli space whether we speak of polarisations of
type $(1,3)$ or $(2,6)$. We denote by $\bAl$ a toroidal
compactification of $\Al$. In this paper it will not matter which
toroidal compactification we take but in [HNS] we shall be more
specific.

If $(A,H,\alpha)$ is a general $(1,3)$-polarised abelian surface with
a level-$2$ structure~$\alpha$ then (cf [BN],~6.1) the blown-up Kummer
surface $\BKm A$ has an embedding into $\PP^3$ whose image is a
Heisenberg-invariant quartic surface.

\thm 1.2 Suppose $(A,H,\alpha)\in\Al$. Let $\sigma:\tilde A\to A$ be
the blow-up of $A$ in the sixteen 2-torsion points and let $f:\tilde
A\to\BKm A = \tilde A/\tilde\iota$ be the quotient by the involution
on $\tilde A$ induced by $\iota$ on $A$. Take $\cL=\cL_{2H}$ as in
Lemma~1.1 and let $\tilde\cL=\sigma^*\cL$. Then the linear system
$|\cL|^-$ of $\tilde\iota$-anti-invariant sections induces a rational
map $\pL:\BKm A\rat\PP^3$: as long as $(A,H)$ is not a product, $\pL$
is a morphism and the image is a quartic surface containing a
line. Moreover, $\alpha$ induces an $\Hl$-action on $\BKm A$ and
$\tilde\cL$ and $\pL$ is equivariant for these actions.

\pf This is merely an assemblage of known results, restated here in a
form convenient for us. The existence of the map $\pL$ may be found in
[Na], [Ni1] and [Ba]. In all these places it is also shown that
the image is a quartic surface. The choice of a symmetric line bundle
$\cM$ such that $\cM^2\imic\cL$ determines a symmetric divisor $M_0$
on $A$ given by the vanishing of an anti-invariant section of $\cM$.
The image of $\sigma^*M_0$ in $\PP^3$ is a line ([BN], p.194) and the
$\Hl$-equivariance is also proved in [BN].~\qed

\corollary 1.3 There is a dominant rational map $\psi:\bAl\rat N$
which is a morphism on $\Al$.

\pf The map $\psi$ is defined by $\psi(A,H,\alpha)=\pL(\BKm A)$,
which, by Theorem~1.2 and [BN], Theorem~8.1, gives a point of $N$ as
long as the 16 lines coming from the 16 choices for $\cM$ are skew.
This is true for general $A$ so $\psi$ extends to an open subset of
$\bAl$ and also the general point of any of the Humbert surfaces that
parametrise product abelian surfaces. We shall see below (Theorem 3.2)
that these Humbert surfaces are contracted to points so $\psi$ extends
to every point of the Humbert surfaces and hence to the whole of
$\Al$. The closure of the image of $\psi$ is an irreducible subvariety
of $N$ and has dimension 3, so $\psi$ is dominant.~\qed

In [HNS] we describe an extension of $\psi$ to part of the boundary.

\startsection 2 Moduli

In this section we shall describe the relationships between the
projective varieties such as $N$ and $M$ on the one hand and moduli
spaces for abelian and Kummer surfaces on the other. We begin by
restating a main result from~[BN].

\prop 2.1 There is a double cover $\tilde N=M/\Hl$ of $N$ such that
$\tilde N$ is birationally equivalent to $\bAl$.~\qed

If we identify $(\sum x_i=0)$ with the Pl\"ucker quadric of lines
in~$\PP^3$ then $M$ parametrises those lines that lie on some
$\Hl$-invariant quartic surface $X\contin\PP^3$. For a general $\ell\in
M$ this $X$ is unique. The action of $\Hl$ on $M$ is described
in~[BN]: every element of $\Hl$ changes the signs of an even number of
the coordinates~$x_i$. Of course $-1\in\Hl$ acts trivially in
projective space, so $\Hl$ acts on $M$ via the quotient and the
morphism $M\to\tilde N$ is of degree~$16$.

The squaring map $M\to N$ of degree~$32$ may be interpreted as the
quotient map under the action of the group generated by $\Hl$ and
$\epsilon:(x_0:x_1:x_2:x_3:x_4:x_5)\mapsto(-x_0:x_1:x_2:x_3:x_4:x_5)$.
It factors through $M\to\tilde N$, the remaining part being the double
cover $\nu:\tilde N\to N$ given by taking the double cover
$\nu:\tilde\PP^5\to\PP^5$ branched along the hyperplane sections, so
$\tilde N=\nu^{-1}(N)$. Thus $\nu$ is induced by the squaring map: as
is pointed out in [BN], the spaces $\bAl$ and $\tilde N$ are
birationally equivalent but the generically $2$-to-$1$ rational map
$\psi:\bAl\rat N$ is completely different from~$\nu$.

We now restate the theorem from [BN] giving moduli descriptions
of~$M$ and~$\tilde N$.

\thm 2.2 $\tilde N$ is birationally equivalent to a compactification
of the moduli space $\Al$ of abelian surfaces~$A$ with a
polarisation~$H$ of type $(1,3)$ and a level-$2$
structure~$\alpha$. $M$~is birationally equivalent to a
compactification of the moduli space of abelian surfaces with a
symmetric bundle~$\cL$ of type~$(1,3)$ and a level-$2$ structure.~\qed

\noindent Remarks. i) The result proved in [BN] is rather more
precise, specifying open sets $M^s$ and $M^s/\Hl$ in $M$ and $\tilde
N$ and open sets in the moduli spaces where the birational
equivalences are isomorphisms.

ii) In view of Lemma~1.1, an alternative birational way of describing
$\tilde N$ is as the moduli space of abelian surfaces with a totally
symmetric line bundle $\cM$ of type~$(2,6)$ and a level-$2$
structure.

iii) There is no need to speak of the moduli of $(2,6)$-polarised
abelian surfaces. The version of this theorem in [Hun] (Theorem~3.4.16
and Corollary~3.4.17) is wrong, because $\Al$ and $\cA_{2,6}(2)$
are the same (and so are $\Gamma_{1,3}$ and $\Gamma_{2,6}$); however,
Hunt's discussion of the situation is scarcely affected by this slip.

iv) There is no preferred choice of~$\cL$: thus $M$ is a principal
$\Ztwo^{\oplus 4}$-bundle over~$\tilde N$. Over a point of $M^s/\Hl$
which is the moduli point of $(A,H,\alpha)$ the fibre may be thought
of as the principal ${}_2A$-space of all (sixteen) symmetric line
bundles of Chern class~$H$, acted on by translation.

Strictly speaking the procedure used in [BN] to prove Theorem~2.2
constructs a family of abelian torsors, not abelian surfaces,
over~$M^s$, since Nikulin's construction does not select a
distinguished origin but only picks out the sixteen fixed points
of~$-1$. But this is good enough because one can replace the family by
its double dual.

The moduli space
$\cA_{1,3}$ of $(1,3)$-polarised abelian surfaces is $\HH_2/\Gpara$ where
$$
\HH_2=\left\{Z=\pmatrix{\tau_1&\tau_2\cr \tau_2&\tau_3\cr}\in
M_{2\times 2}(\CC)\mid Z={}^tZ,\quad\im Z>0\right\}
$$
and $\Gpara$ is the paramodular group in the sense of [GH], that is
$$
\Gpara=\left\{\gamma\in\Sp(4,\QQ)\mid \gamma \in \pmatrix{
\ZZ&\ZZ&\ZZ&3\ZZ\cr
3\ZZ&\ZZ&3\ZZ&3\ZZ\cr
\ZZ&\ZZ&\ZZ&3\ZZ\cr
\ZZ&{{1}\over{3}}\ZZ&\ZZ&\ZZ\cr
}\right\}.
$$
$\Gpara$ acts on $\HH_2$ by fractional linear transformations,
$\pmatrix{A&B\cr C&D\cr}:Z\mapsto (AZ+B)(CZ+D)^{-1}$. 

In the case of polarisation of type $(1,3)$ (or type $(1,t)$ in fact)
it makes sense to speak of a dual polarised variety: according to [GH]
one has a map $\Phi(3):\cA_{1,3}\to\cA_{1,3}$ which sends
$(A,H)$ to $(\hat A,\hat H)$ where $\hat A=\Pic^0 A$ and $\hat H$ is
of type~$(1,3)$. It is induced ([GH], Proposition~1.6) by the element
of order~$2$
$$
V_3=\pmatrix{0&{\sqrt 3}^{-1}&0&0\cr
             \sqrt 3&0&0&0\cr
             0&0&0&\sqrt 3 \cr
             0&0&{\sqrt 3}^{-1}&0\cr}
\in\Sp(4,\RR)
$$
also acting on $\HH_2$ by a fractional linear transformation.

$\Al$ is the quotient $\HH_2/\Gl$, where $\Gl<\Sp(4,\ZZ)$ is given by
$$
\Gl=\left\{\gamma\in\Sp(4,\QQ)\mid \gamma-I \in \pmatrix{
2\ZZ&2\ZZ&2\ZZ&6\ZZ\cr
6\ZZ&2\ZZ&6\ZZ&6\ZZ\cr
2\ZZ&2\ZZ&2\ZZ&6\ZZ\cr
2\ZZ&{{2}\over{3}}\ZZ&2\ZZ&2\ZZ\cr
}\right\}.
$$
$V_3^{-1}\Gl V_3=\Gl$ so the involution $\Phi(3)$ induces an
involution on $\Al$ with the same properties: we shall call this
$\Phi(3)$ also. 

There is an alternative moduli description of $\tilde N$ as a space of
torsors. Because we start with an $H_{2,2}$ orbit of lines, rather
than a particular one, this is what the geometry really produces.

\prop 2.3 $\tilde N$ is birationally equivalent to a compactification
of the moduli space of symmetric abelian torsors~$Y$ of dimension~$2$ with a
symmetric bundle~$\cL$ of type~$(1,3)$ and a level-$2$ structure~$\alpha$.

\pf Given a triple $(A,H,\alpha)\in\Al$, where $A$ is an abelian variety we
associate to it $(\hat A,\cL,\hat\alpha)$, where $\hat A$ is the torus
obtained by forgetting the origin of~$A$, and $\cL$ is a symmetric
line bundle representing~$H$. If we choose another symmetric line
bundle~$\cL'$ with $c_1(\cL')=H$, then there is a 2-torsion
point~$\bdx$ with $t_\bdx^*\cL=\cL'$. This translation defines an
automorphism of the torsor and hence this map is well defined up to
isomorphism. The level-$2$ structure $\hat\alpha$ is simply a
level-$2$ structure on the underlying abelian surface. Inversely,
starting with $(Y,\cL,\hat\alpha)$ with $Y$ a symmetric torsor, we
make $Y$ into an abelian surface $A$ by choosing one of the sixteen
fixed points of the involution~$-1$, say $\bdx$, as the origin. We
take $H=c_1(\cL)$ and $\alpha=\hat\alpha$. We must check that this
definition does not depend on the choice of $\bdx$. If we choose a
different $2$-torsion point $\bdx'$ instead of $\bdx$ we get an
abelian surface $A'$, and the translation $t_{\bdx'-\bdx}$ defines an
isomorphism from $A$ to~$A'$. It maps~$\cL$ to a different line
bundle~$\cL'$, but $c_1(\cL)=c_1(\cL')$ so again at the level of
isomorphism classes the map is well defined and the two maps are
obviously inverse to each other.~\qed

The above argument works for any polarisation (and indeed any
dimension). In our case it is also possible to use $\Phi(3)$ and the
totally symmetric bundle of class~$2H$ to induce an isomorphism.

A symmetric bundle of type $(1,3)$ even determines a unique divisor on
$Y$, because the space of anti-invariant sections under the action
of~$-1$ is $1$-dimensional. So $Y$ is very like a Jacobian but the
``theta divisor'' gives a non-principal polarisation.

The relevance of $\Phi(3)$ is that it gives rise to~$\psi$.

\prop 2.4 There exists $g\in\Gpara/\Gl$ such that if $x_1$,
$x_2\in\Al$ and $\psi(x_1)=\psi(x_2)$ then $\Phi(3)(gx_1)=x_2$.

\pf Suppose $x_i=(A_i,H_i,\alpha_i)$ are such that $\psi(x_1)=\psi(x_2)$. For
general $x_i$ we know that $\Km A_1\imic\Km A_2$, since $\pL$ is an
isomorphism onto its image. By Theorem~1.5 of~[GH], we also know that,
for general $A_i$, if $\Km A_1\imic\Km A_2$ then $(A_2,H_2)=(\hat
A_1,\hat H_1)$. So the map $\Al\to\Al$ induced by $\psi$ agrees with
$V_3$ up to a possible change of the level structure, that is, up to
the action of~$\Gpara/\Gl$.~\qed
 
The element $gV_3$ is an involution modulo $\Gamma_{1,3}(2)$. There
are, in principle, two possibilities. One is that $gV_3$ itself is an
involution. (This happens if $g=1$, but also for some other~$g$.) By
[GH], Corollary~3.9, there is exactly one involution in
$\Gpara V_3$ up to conjugation with $\Gpara$. The fixed locus of
such an involution is a Humbert surface of discriminant~$12$ in
$\Al$. For a general element of this Humbert surface the linear system
$\phi_{|\tilde\cL|^-}$ is an embeddding. But then this would
contradict the statement~8.4 of~[BN], which says that $\psi$ is
unbranched on $M^s/\Hl$. In particular we see that $g \neq 1$. The
other possibility is that $(gV_3)^2 \in \Gamma_{1,3}(2)$, but
$(gV_3)^2 \neq 1$. It follows, again by [GH], Corollary~3.9 that the
map~$\psi$ is not branched along a Humbert surfaces of
discriminant~$4$. On the other hand the branch locus of the map~$\nu$
is the image of the locus of bielliptic abelian surfaces (see
below). These surfaces are parametrized by a Humbert surface of
discriminant~$4$. Hence we can see directly from moduli that the two
maps $\psi$ and~$\nu$ do not agree.

\startsection 3 Branching and special abelian surfaces

In this section we study the branch locus of $\nu$ below $\Al$ and the
corresponding abelian surfaces. Recall that we define
$\psi:\Al\to|\cO_{\PP^3}(4)|$ by
$\psi(A,H,\alpha)=X=\im\pL\contin\PP^3$. By continuity
$\psi(A,H,\alpha)\in N$, even if $X$ is singular.

The image of the map $\pL$ was studied by Bauer [Ba]. There is a
classical study of quartic surfaces in $\PP^3$ with singularities by
Jessop [J]. We begin, though, with the double points of $N$. These are
in the S-planes but, unlike other points in the S-planes, they can
arise from smooth abelian surfaces. According to the main theorem of
[Ba] the map $\pL:\BKm A\to\PP^3$ is usually an isomorphism onto its
image (this is also shown in [BN] and [Na]) but there are exceptional
cases, numbered (III), (IV) and (V) in \S{}5 of~[Ba].

\noindent{\bf III.} The $(1,3)$-polarisation $H$ is given by
$\cO_A(G+E)$ with $E$ elliptic and $G$ an irreducible genus~$2$ curve
with $G.E=2$. In this case $\pL$ is birational and the image $X$ has
four nodes.

\noindent{\bf IV.} $H$ is given by $\cO_A(E_1+E_2+E_3)$, $E_i$
elliptic, $E_i.E_j=1-\delta_{ij}$. This is a degenerate case of (III)
and the image $X$ has twelve nodes. $\pL$ is still birational.

\noindent{\bf V.} $A=E_1\cross E_2$ and $H=c_1(\cO_A(E_1+3E_2))$,
$E_i$ elliptic. The image $X$ is smooth but $\pL$ is $2$-to-$1$.

\lemma 3.1 If a group $G$ acts transitively on a set $X\ni x_0$ and
$G'$ is a normal subgroup of $G$ of finite index, and if $S$ and $S'$
are the stabilisers of $x_0$ in $G$ and $G'$ respectively, then the
number of $G'$-orbits in $X$ is equal to $[G:G']/[S:S']$.

\pf This is an elementary calculation (cf. Lemma~2.2 of~[HW]).~\qed

\thm 3.2 The nodes in $N$ correspond to points of $\Al$ where
$A=E_1\cross E_2$ and $H=c_1(\cO_A(E_1+3E_2))$. These points form a
surface with twenty irreducible components, which are
contracted by $\Al\to N$ to the 10 nodes of $N$.

\pf The 10 nodes of $N$ are the $\SS_6$-translates of
$(1:1:1:-1:-1:-1)\in \PP^4$ and they correspond to nonreduced quartic
surfaces, namely the squares of the fundamental quadrics. This is
checked for the point $(-1:1:1:1:-1:-1)$ (which is enough) in [BN],
p.190. As there are no other Heisenberg-invariant quadrics, any
nonreduced $\Hl$-invariant quartic in $\PP^4$ is the square of a
fundamental quadric and corresponds to a node of~$N$.

Suppose $A=E_1\cross E_2$ and $H=c_1\big(\cO_{E_1}(3) \sqtimes
\cO_{E_2}(1)\big)$ is a product polarisation of type $(1,3)$ on $A$. Then
we are in case~(V) of [Ba] and therefore by the main theorem of [Ba]
the rational map $\pL:\tilde A\rat \PP^3$ is not defined on certain
base curves of $|\tilde{\cL}|^-$, namely the symmetric translates of
$E_2$. However, away from these base curves $\pL$ coincides with the
morphism $\tilde A\to\PP^3$ coming from the polarisation
$c_1\big(\cO_{E_1}(2)\sqtimes\cO_{E_2}(2)\big)$ of type $(2,2)$, which
is of degree 2. So the closure of the image of $\pL$ in $\PP^3$ is a
(double) quadric: if we make a choice of level-$2$ structure so as to
fix an $\Hl$-action the image will then be $\Hl$-invariant so it must
be one of the ten fundamental quadrics. So the image of the locus
$$
\big\{(A,H,\alpha)\in\Al\mid(A,H)\imic(E_1\cross E_2,
c_1(\cO_{E_1}(3)\sqtimes\cO_{E_2}(1)))\big\},
$$
of product surfaces with product polarisation, in $N$ is contained in
the 10 nodes. Since the nodes are permuted transitively by the action
of $\SS_6$, the image must be all the nodes.

Conversely, if $\psi(A,H,\alpha)$ is one of the nodes then
$(A,H,\alpha)$ must be a product since (by [Ba]) in other cases the
map $\pL:\BKm(A)\to\PP^3$ is birational onto its image, which is
therefore irreducible and birationally a K3 surface.

To prove the rest of the theorem we need to count the number
of irreducible components of the space
$$
\big\{(A,H,\alpha)\mid (A,H)\hbox{ a product, }H\hbox{ of type
}(1,3),\ \alpha\hbox{ a level-$2$ structure}\big\}
$$
in $\Al$. We use the same method as [HW], making use of the
description of $\Al$ as a Siegel modular variety. 

Consider the surface $\hat{\cH}_1=\{Z\in\HH_2\mid \tau_2=0\}$ and its
image $\cH_1$ in $\cA_{1,3}$. It is shown in [HW] that $\cH_1$
parametrises abelian surfaces of type $(1,3)$ which split as polarised
abelian surfaces, that is, product surfaces with the product
polarisation. $\cH_1$ is irreducible, because $\hat\cH_1$ is. We want
to know the number of components of $\pi^*(\cH_1)$, where
$\pi:\Al\to\cA_{1,3}$ is the quotient map by the action of
$\Gpara/\Gl\imic\Sp(4,\ZZ_2)$ (which is $\SS_6$ as an abstract
group).

We look at the stabilisers
$S_1=\{\gamma\in\Gpara\mid\gamma(\Hat\cH_1)=\Hat\cH_1\}$
and $S_1(2)=\{\gamma\in\Gl\mid\gamma(\Hat\cH_1)=\Hat\cH_1\}$. As in
[HW] one has
$$
S_1=\left\{\pmatrix{
a&0&b&0\cr
0&a'&0&b'\cr
c&0&d&0\cr
0&c'&0&d'\cr
}\mid\pmatrix{
a&b\cr c&d&\cr},
\pmatrix{a'&b'\cr c'&d'\cr}
\in\SL(2,\ZZ)\right\}
$$
so $S_1\imic\SL(2,\ZZ)\cross\SL(2,\ZZ)$, and an element of $S_1$ is in
$\Gl$, and hence in $S_1(2)$, if and only if $\pmatrix{a&b\cr c&d\cr}$
and $\pmatrix{a'&b'\cr c'&d'\cr}$ are both congruent to the identity
mod~2. So $S_1/S_1(2)\imic\SL(2,\ZZ_2)\cross\SL(2,\ZZ_2)$.

Now we can conclude the proof of Theorem 3.2 by applying Lemma 3.1. We
take $G=\Gpara$, $G'=\Gl$, and $X=\{\gamma\Hat\cH_1\mid
\gamma\in\Gpara\}$, the set of $\Gpara$-translates of $\Hat\cH_1$
or of preimages of $\cH_1$ in $\HH_2$. The set $X/G'$ is precisely the
set of preimages of $\cH_1$ in $\Al$, and by Lemma 3.1 it has
$$
[\Gpara:\Gl]/[S_1:S_1(2)]=|\Sp(4,\ZZ_2)|/|\SL(2,\ZZ_2)\cross\SL(2,\ZZ)|
=720/36=20
$$
elements. This is therefore the number of irreducible surfaces in
$\Al$ corresponding to product surfaces.~\qed

The nodes lie in the S-planes in $N$. We have said that in general
points in the S-planes will correspond to degenerate abelian
surfaces, but here the closure of an Humbert surface in $\bAl$ is
contracted by $\psi$ to a point which we may also think of as coming
from the boundary $\bAl\setminus\Al$. Of course other components in the
boundary may also be contracted to the nodes (this will depend on our
choice of toroidal compactification), but only the nodes come from
both degenerate and nondegenerate abelian surfaces.

\thm 3.3 If $\bdu\in N$ is a point in some S-plane and there is an
abelian surface $(A,H,\alpha)\in\Al$ such that
$\psi(A,H,\alpha)=\bdu$, then $\bdu$ is a node and $(A,H,\alpha)$ is a
product.

\pf If $(A,H,\alpha)$ is not a product then $\pL:\BKm(A)\rat X$ is
birational onto its image, which is therefore an irreducible variety
which is birationally a K3 surface. But for any $\bdu$ on a S-plane
the corresponding quartic surface contains a double line (actually a
pair of double lines). Following [J], Chapter~VI, \S{}77, we observe
that a plane containing this line cuts out a conic on~$X$. So $X$ is
covered by a family of rational curves and is therefore ruled.~\qed

Now we turn to the V-planes $F_{ij}$, defined by $u_i=u_j=0$,
which form the branch locus of $\nu:\tilde N \to N$. There are fifteen of
these planes, equivalent under the $\SS_6$-action on $\PP^5$: we may
choose one of them to work with and we use~$F_{45}$. $F_{45}$
determines the family
$$
\cF_{AE}=\{B(z_0^2z_1^2+z_2^2z_3^2)+C(z_0^2z_2^2+z_1^2z_3^2)
+D(z_0^2z_3^2+z_1^2z_2^2) = 0\}
$$
which is studied in [Ni2].

\thm 3.4 Let $\bdu\in N$ be a point of a V-plane, not lying on an
S-plane nor on one of the D-lines $\{u_i=u_j=u_k=0\}$. Then there
exists an abelian variety $A$ with a polarisation $H$ of type $(1,3)$
and a level-$2$ structure~$\alpha$, such that $\psi(A,H,\alpha)=\bdu$
and $(A,H)$ is a bielliptic abelian surface in the sense of [HW]. In
particular $A$ is isomorphic to $E_1\cross E_2/\ZZ_2\cross\ZZ_2$ for
some elliptic curves $E_1$, $E_2$ and some action of
$\ZZ_2\cross\ZZ_2$.

\pf The conditions on $\bdu$ imply that the corresponding quartic
surface $X_\bdu=X$ (given by $B=\hf(-u_0-u_1+u_2+u_3)$, etc; see
[BN], [Ni2]) has four simple nodes and no other singularities. As in
[Ni2] we can use the theorem of Nikulin [N] to construct a diagram
$$
A{\buildrel \sigma \over \longleftarrow}\tilde
A {\buildrel f \over \To} \tilde X {\buildrel \beta \over \To} X
$$
where $\beta$ is the blow-up of the four nodes of~$X$ and $f:\tilde A
\to \tilde X$ is the double cover branched along a set ${\bf L}_0$ of
sixteen disjoint smooth rational curves. ${\bf L}_0$ is the pullback
to $\tilde X$ of one of the two $\Hl$-orbits of lines in~$X$, which
become disjoint after the blow-up~$\beta$. We consider the curve
$E_0=\sigma_*f^*\big(\beta^{-1}(p_0)\big)$, where $p_0$ is one of the
nodes, and $\Theta=\sigma_*f^*\beta^*L'$, where $L'$ is another line
on~$X$. Then ([Ni2], Lemma~6-3) $g(E_0)=1$ and $g(\Theta)=2$.

For reasons explained below we now diverge from the line of argument
given in [Ni2]. We have $E_0.\Theta=2$ and there is an exact sequence
$$
0\To E_0 \To A \To E'_0 \To 0
$$
for some elliptic curve~$E'_0$. So the map $\pi:A\to E'_0$ induces a double
cover $\Theta\to E'_0$. By Torelli, $(A,\cO_A(\Theta))=(\Jac\Theta,\Theta)$
so $A$ is the Jacobian of a bielliptic genus~2 curve, i.e., a principally
polarised bielliptic abelian surface. By [HW], Proposition~4.1, it follows
that $A=E_1\cross E_2/\ZZ_2\cross\ZZ_2$ for suitable elliptic curves $E_1$,
$E_2$.

The polarisation which we are really interested in, however, is
$H=E_0+\Theta$, which is of type~$(1,3)$. The bielliptic involution
$j:\Theta\to\Theta$ induces $j_A:(A,\Theta)\to(A,\Theta)$ and
$\pi:A\to E'_0$ factors through $j_A$. Therefore $j_A$ preserves
$E_0$, which is a fibre of $\pi$, so $j_A$ induces an involution of
$(A,H)$ as well. Hence (cf. [HW]) $(A,H)$ either splits as a
polarised abelian surface or else is a bielliptic $(1,3)$-polarised
abelian surface. The first case can be excluded, because we already
know that the linear system $|2H|^-$ induces a map of $\BKm A$ which
is birational onto its image. (Looking carefully at the proof of
Proposition~4.4 of [HW] one can show that every bielliptic
$(1,3)$-polarised abelian surface contains a bielliptic genus~$2$
curve, and conversely.) The level-$2$ structure on $(A,H)$ is simply
induced by the $\Hl$-action on~$X$.~\qed

\noindent{\it Remark.}\/ In [Ni2] the second elliptic curve and the
$\ZZ_2\cross \ZZ_2$-action are constructed more directly, but under the
additional assumption that $\rho=\rk\NS(A) =2$. Also, there is a gap in the
proof of Lemma~6-3(2) of [Ni2]. Nevertheless, we could have used this
method. The restriction $\rho=2$ is harmless to us because it is true for
almost all bielliptic surfaces, and the gap can easily be filled. We take
the opportunity to do this.

In $\NS(A)$ we have ([Ni2], p.333) inequivalent elliptic curves $E$, $E'$
and a genus~2 curve $\Theta'$ with $E.\Theta'=2$. So in $\NS(A)\tens\QQ$ we
have $E'=n_0E+n_1\Theta'$ and therefore $n_1=-2n_0\not=0$. It does not
follow at once, as asserted in [Ni2], that $n_0=\pm 1$, because we do not
know that $n_0\in\ZZ$. The elliptic curves, however, are not divisible in
$NS(A)$, because topologically $A=E\times T$ for some real subtorus
$T\subset A$. So the denominator of $n_0$ is at most~$2$,
otherwise $E'$ is divisible, and if $n_0\in\ZZ$ then $n_0=\pm 1$. So
suppose $n_0=k+{{1}\over{2}}$. Then $2E'=(2k+1)E-(4k+2)\Theta'$ so
$E=2\big(E'+(2k+1)\Theta'-kE\big)$ is divisible. So $n_0$ and $n_1$ are in
fact integers and the argument in~[Ni2] can be used.

\corollary 3.5 The closures of the following three loci in $\bAl$ coincide:
\hfil\break\hskip 1 cm (i) the locus of $(A,H,\alpha)$ such that $(A,H)$ is a
$(1,3)$-polarised bielliptic abelian surface;
\hfil\break\hskip 1 cm (ii) the locus of $(A,H,\alpha)$ such that $\pL(\BKm
A)=X_\bdu$ for some $\bdu$ in a V-plane, not in a S-plane or
on a D-line;
\hfil\break\hskip 1 cm (iii) the locus of $(A,H,\alpha)$ such that
$H=c_1(\cO_A(E+\Theta))$ with $E.\Theta=2$ and $\Theta\contin A$ a
symmetric irreducible curve of genus~2.
\hfil\break Moreover, the loci in (ii) and (iii) coincide precisely.

\pf We have shown that (ii) implies (i) in Theorem~3.4. By the main
theorem of [Ba], the hypotheses of (iii) imply that $|\tilde\cL|^-$
contracts the images in $\BKm A$ of the four symmetric translates of
$E$ and no more, so that $\pL(\BKm A)$ is an $\Hl$-invariant quartic
surface with 4 nodes. These are precisely those parametrised by the
points of the V-planes with the exceptions given in (ii), so (ii) and
(iii) are the same. Finally, the locus of bielliptic surfaces without
level structure in $\cA_{1,3}$ and the locus of Jacobians of
bielliptic genus~2 curves in the moduli space of principally polarised
abelian surfaces are both irreducible surfaces, as is shown in
[HW]. The locus of such Jacobians can also be thought of as a
subvariety of $\cA_{1,3}$, by replacing the principal polarisation
$\Theta$ by $E_0+\Theta$ as above: there is a unique way to do
this. The two irreducible surfaces in $\cA_{1,3}$ must therefore
coincide, so the closures of the loci in (i) and (iii) are the
same. In fact the loci in (i) and (iii) coincide precisely of we drop
the irreducibility condition from (iii), as one can deduce from the
proof of Proposition~1.4 of~[HW].~\qed

Now we know that the V-planes in $N$ are precisely the image under
$\psi$ of the surface in $\Al$ parametrising bielliptic surfaces
of type~$(1,3)$; that is, of $\pi^{-1}(\cH_2)$, where
$\pi:\Al\to\cA_{1,3}$ is the quotient map by $\Sp(4,\ZZ_2)$ as before
and $\cH_2\contin\cA_{1,3}$ parametrises bielliptic abelian
surfaces. We can count the irreducible components of $\pi^{-1}(\cH_2)$
just as we did for $\pi^{-1}(\cH_1)$.

\thm 3.6 The surface in $\Al$ parametrising triples $(A,H,\alpha)$ for
which $(A,H)$ is bielliptic has 15 irreducible components.

\pf We put $\hat\cH_2=\{Z\in\HH_2\mid 3\tau_1=2\tau_2\}$ as in [HW].
We know, also from [HW], that the image of $\hat\cH_2$ in $\cA_{1,3}$
is $\cH_2$, so according to Lemma~3.1 and the same reasoning as in
Theorem~3.2 the number of components of $\pi^{-1}(\cH_2)$ is
$$
[\Gpara:\Gl]/[S_2:S'_2]=720/[S_2:S_2(2)]
$$
where $S_2$ and $S'_2$ are the stabilisers of $\hat\cH_2$ in
$\Gpara$ and $\Gl$ respectively. From [HW] we have
$$
\eqalign{
S_2&=\Bigg\{
\pmatrix{
a&0&2b&3b\cr
3(a-a')/2&a'&3b&3(3b+b')/2\cr
(3c'+c)/2&-c'&d&3(d-d')/2\cr
-c'&2c'/3&0&d'\cr
 }
\cr
&\qquad\qquad\Bigg|\> M=\pmatrix{a&b\cr c&d\cr}, M'=\pmatrix{a'&b'\cr
c'&d'\cr}\in\SL(2,\ZZ), M\equiv M'\hbox{ mod }2\Bigg\}\cr
   &\imic\big\{(M,M')\in\SL(2,\ZZ)\cross\SL(2,\ZZ)\mid M\equiv
M'\hbox{ mod }2\big\}.\cr}
$$
Such a matrix will be in $\Gl$ if and only if $M\equiv M'\equiv I$
mod~2 and $M\equiv M'$ mod~4, because we need $a$, $a'$, $d$, $d'$
odd, $b$, $c'$ even (so $b'$ and $c$ also even) and
$3(a-a')/2$, $3(3b+b')/2$, $(3c'+c)/2$ and $3(d-d')/2$ all even.

Now consider the tower of groups
$$
\eqalign{
\tilde S_2=\{(M,M')&\mid M\equiv M'\hbox{ mod }2\}\cr
&\bigcup\cr
\tilde S''_2=\{(M,M')&\mid M\equiv M'\equiv I\hbox{ mod }2\}\cr
&\bigcup\cr
\tilde S'_2=\{(M,M')&\mid M\equiv M'\equiv I\hbox{ mod 2, }M\equiv
M'\hbox{ mod }4\}.\cr}
$$
Then the map
$$
(M,M')\mapstochar\To
\pmatrix{
a&0&2b&3b\cr
3(a-a')/2&a'&3b&3(3b+b')/2\cr
(3c'+c)/2&-c'&d&3(d-d')/2\cr
-c'&2c'/3&0&d'\cr}
$$
induces isomorphisms between $S_2$ and $\tilde S_2$ and $S'_2$ and
$\tilde S'_2$, so $[S_2:S'_2]=[\tilde S_2:\tilde S'_2]$. Now $\tilde
S''_2$ is a normal subgroup of $\tilde S_2$ and the quotient is
isomorphic to $\SL(2,\ZZ_2)$, of order~6, and the cosets of $\tilde
S'_2$ in $\tilde S''_2$ are determined by the residue class of $M$
mod~4. There are eight of these: the off-diagonal elements may be 0~or
2~mod~4, and the diagonal ones are either both~$1$ or both~$-1$. So
$[\tilde S''_2:\tilde S'_2]=8$, so $[\tilde S_2:\tilde S'_2]=[\tilde
S_2:\tilde S''_2][\tilde S''_2:\tilde S'_2]=48$. So the number of
irreducible components is $720/48=15$.~\qed

Finally, we consider the D-lines $L_{ijk}=\{u_i=u_j=u_k=0\}\contin N$
and examine the abelian surfaces that live over them. We work with the
line $L_{045}$. On this line the quintic surface $X$ has the
equation
$$
C(z_0^2-z_3^2)(z_1^2-z_2^2)+D(z_0^2-z_2^2)(z_3^2-z_1^2)=0.
$$
These are the surfaces known as desmic\footnote {${}^2$}{{\smallfont Greek
$\scriptstyle{\delta\muskip-0.5mu\epsilon\muskip-0.5mu\sigma\muskip-0.5mu
\mu\muskip-0.5mu{\skew2\acute o}\muskip-0.5mu\varsigma}$,
a band or tie, from
$\scriptstyle{\delta\muskip-0.5mu\skew2\acute\epsilon\muskip-0.5mu\omega}$,
to bind.}} surfaces in [Hud] and~[J].

The surface $\SCD$ given by the equation above is singular at twelve
points and contains 32 lines. Sixteen of these lines are common to all
the surfaces in the pencil. More precisely we have the following
result.

\thm 3.7 {\rm ([Hud], [J])} If $CD\not=0$ then $\SCD$ has ordinary
double points at the four poles $P_0,\ldots,P_3$ ($=(1:0:0:0)$, etc)
and at the eight points $(\pm 1:\pm 1:\pm 1:\pm 1)$ forming the
vertices of a cube, and no other singularities. A vertex is said to be
even or odd according to whether the number of coordinates equal to
$-1$ is even or odd. The sixteen lines $L_1,\ldots,L_{16}$ which are
the edges and main diagonals of the cube lie in $\SCD$: each passes
through exactly one pole, one even point and one odd point.

The configuration of lines and points is the well-known Reye
configuration: see also in this context [Ni2], \S\S{}6--7.

We can resolve the singularities of $\SCD$ by blowing up the nodes.
Let $\beta:\tS\to\SCD$ be this blow up and let $\tilde L_j$ be the
proper transform of $L_j$ in $\tS$. Denote by $C_i^0$, $C_i^+$ and
$C_i^-$ the exceptional curve in $\tS$ coming from the $i$th pole,
even and odd vertex respectively.

\lemma 3.8 $\tS$ is a smooth K3 surface.

\pf $K_\SCD=0$ by adjunction, and blowing up nodes does not change this.~\qed

The sixteen rational curves $\tilde L_j$ are smooth and disjoint, so
we can apply Nikulin's construction from [N] as before, obtaining
$$
A_{(C:D)}{\buildrel \sigma \over \longleftarrow}\tilde
A_{(C:D)}{\buildrel f \over \To} \tS {\buildrel \beta \over \To}\SCD
$$
where $f$ is a double cover branched along $\tilde L_j$ only, $\sigma$
is the blow-down of the sixteen rational curves $f^{-1}(\tilde L_j)$
and $A_{(C:D)}$ is (after choosing an origin) an abelian surface.
Furthermore, $\tS=\BKm A_{(C:D)}$ and $f$ is the Kummer map.

\thm 3.9 $A_{(C:D)}\imic E\cross E$, where $E$ is the elliptic curve
whose $j$-invariant is $j(C/D)=2^8(C^2-CD+D^2)^3/C^2D^2(C-D)^2$.

\pf Consider the curves $\tilde E_i^0=f^{-1}(\tilde C_i^0)$, etc. Each
of these is an elliptic curve since $f:\tilde E_i^0\to\tilde C_i^0$ is
a double cover branched at the four points corresponding to the four
lines $L_j$ passing through $P_i$, and similarly for the odd and even
vertices. These twelve elliptic curves are disjoint in $\tilde
A_{(C:D)}$ but their images $E_i^0$, $E_i^\pm$ in $A_{(C:D)}$ have
intersection numbers $E_i^0E_j^0=0$,
$E_i^0E_j^+=E_j^+E_k^-=E_k^-E_i^0=1$. Indeed $E_i^0$, $E_j^+$ and
$E_k^-$ are concurrent.

There is an exact sequence of abelian varieties
$$
0\To E_0^0\To A_{(C:D)}\To E_0^{0*}\To 0
$$
where $E_0^{0*}$ is an elliptic curve. Since $E_0^0E_j^+=1$ it follows
that $E_j^+\imic E_0^{0*}$ and similarly $E_k^-\imic E_0^{0*}$ for all
$j$, $k$. Similarly $E_i^0\imic E_0^{+*}\imic E_0^-$ for all~$i$, so
all the elliptic curves $E_i^0,E_j^+,E_k^-,E_i^{0*},E_j^{+*},E_k^{-*}$
are isomorphic to one another. Moreover, $E_0^+$ defines a section of
$A_{(C:D)}\to E_0^{0*}$ so $A_{(C:D)}\imic E_0^0\cross E_0^+\imic
E\cross E$.

To calculate the $j$-invariant of $E$ we need to know the cross-ratio
of the four branch points of $f:\tilde E_0^0\to \tilde C_0^0$. We use
$z_1$, $z_2$ and $z_3$ as affine coordinates in the affine piece
$z_0=1$ of $\PP^3$, so that $P_0$ is the origin of $\AA^3$ and $\SCD$
has the affine equation
$$
C(1-z_3^2)(z_1^2-z_2^2)+D(1-z_2^2)(z_3^2-z_1^2)=0.
$$
The tangent cone to $\SCD$ at $(0,0,0)$ has the equation
$$
C(z_1^2-z_2^2)+D(z_3^2-z_1^2)=0
$$
and this can also be thought of as the equation of $\tilde C_0^0$ in the
$\PP^2$ which is the set of lines in $\AA^3$ through the origin. The branch
points are given by the lines joining $P_0$ to the even points, which
are $(1:\pm 1:\pm 1)$ (the vertices of a square). In the affine piece
of $\PP^2$ given by $z_1=1$ these points are simply $(\pm 1,\pm 1)$
and $\tilde C_0^0$ has the affine equation
$$
C(1-z_2^2)+D(z_3^2-1)=0.
$$
If we take the isomorphism $\tilde C_0^0\to\PP^1$ given by projection
from $(-1,1)$ we find that the four points are mapped to $0$, $1$,
$\infty$ and $C/D$.~\qed

In fact the elliptic curves $E_i^0$, $E_j^+$ and $E_k^-$ are
isomorphic for reasons of projective geometry, because there is an
action of $\SS_6$ that preserves $\SCD$ and permutes the three
tetrahedra parametrised by poles, odd vertices and even vertices.

This is case~(IV) (the diagonal case) of [Ba], \S{}5. The
$(1,3)$-polarisation is given by $\cO(E\cross\{0\}+\{0\}\cross
E+\Delta)$, where $\Delta$ is the diagonal. The result has a long
history: compare [Hun], p.~313, where another proof of part of this
appears, or indeed~[J].

\startsection References

\noindent[Ba] T.~Bauer, {\it Projective images of Kummer surfaces},
Math. Ann. {\bf 299} (1994), 155--170.

\noindent[BN] W.~Barth \& I.~Nieto, {\it Abelian surfaces of type
$(1,3)$ and quartic surfaces with 16 skew lines}, J. Alg. Geom {\bf 3}
(1994), 173--222.

\noindent[GH] V.~Gritsenko \& K.~Hulek, {\it Minimal Siegel modular
threefolds}, Math. Proc. Cam. Phil. Soc. {\bf 123} (1998), 461--485.

\noindent[HNS] K.~Hulek, I.~Nieto \& G.K.~Sankaran, {\it Degenerations
of abelian surfaces of type $(1,3)$ and Kummer surfaces}. To appear in {\it
Proceedings of the Conference on Algebraic Geometry -- Hirzebruch 70}
(P.~Pragacz, M.Szurek \& J. W\'\i{}sniewki, Eds.), AMS Contemporary
Mathematics, 1999.

\noindent[HW] K.~Hulek \& S.~Weintraub, {\it Bielliptic abelian surfaces},
Math. Ann. {\bf 283} (1989), 411--429.

\noindent[Hud] R.W.H.T.~Hudson, {\it Kummer's quartic surface}: Cambridge
University Press, Cambridge 1906 (reissued 1990).

\noindent[Hun] B.~Hunt, {\it The Geometry of some special Arithmetic
Quotients}, Lecture Notes in Mathematics {\bf 1637}: Springer, Berlin 1996.

\noindent[J] C.M.~Jessop, {\it Quartic surfaces with singular points}:
Cambridge University Press, Cambridge 1916.

\noindent[LB] H.~Lange \& Ch.~Birkenhake, {\it Complex abelian varieties}:
Springer, Berlin 1992.

\noindent[Mu] D.~Mumford, {\it On the equations defining abelian
varieties I}, Inv. Math. {\bf 1} (1966), 287--354.

\noindent[Na] I.~Naruki, {\it On smooth quartic embedding of Kummer surfaces},
Proc. Japan Acad. Ser.~A {\bf 67} (1991), 223--225.

\noindent[Ni1] I.~Nieto, {Invariante Quartiken unter der Heisenberg
Gruppe~$T$}, Thesis: Erlangen 1989.

\noindent[Ni2] I.~Nieto, {\it Examples of abelian surfaces with polarization
type $(1,3)$}. In: {\it Algebraic Geometry and Singularities},
(C.~Lopez and N.~Macarro, eds.), Progress in Mathematics {\bf 134},
319--337: Birkhauser, Basel 1996.

\noindent[N] V.~Nikulin, {\it On Kummer surfaces}, Math. USSR. Izv {\bf 9}
(1975), 261--275.

\bigskip
\vbox{\baselineskip=15pt
\font\smalltt=cmtt9
{\settabs3\columns
\+K.~Hulek&I.~Nieto&G.K.~Sankaran\cr
\+Institut f\"ur Mathematik&Cimat, A.C.&Department of\cr
\+Universit\"at Hannover&Callejon de Jalisco S/N&\ Mathematical Sciences,\cr
\+Postfach 6009&Col. Mineral de Valenciana&University of Bath,\cr
\+D 30060 Hannover&36000 Guanajuato, Gto.&Bath BA2 7AY\cr
\+GERMANY&MEXICO&ENGLAND\cr
\+{\smalltt hulek@math.uni-hannover.de}&{\tt nieto@fractal.cimat.mx}&{\tt
gks@maths.bath.ac.uk}\cr}}

\end